\documentclass[11pt,lettersize,reqno,fullpage]{article}

\usepackage{amsmath,amssymb,amsfonts,mathrsfs,verbatim,enumitem,pstricks,amsthm}
\usepackage[all]{xy}
\usepackage{centernot, xr, accents}
\usepackage{hyperref}
 
\title{{\LARGE Enumeration of singular hypersurfaces on  arbitrary complex manifolds}}
\author{Ritwik Mukherjee }
\date{}
\usepackage{hyperref}
\hypersetup{
	colorlinks,
	citecolor=black,
	filecolor=black,
	linkcolor=blue,
	urlcolor=black
}

\def \hf{\hspace*{0.5cm}}                      
\def\bge{\begin{equation}}                
\def\ede{\end{equation}}                
\def\bgd{\begin{displaymath}}         
\def\edd{\end{displaymath}}            
\def\bgee{\begin{equation*}}           
\def\edee{\end{equation*}}           

\def \ni{\noindent}
\def \hf{\hspace*{0.5cm}}                      
\def\lra{\longrightarrow}

\def\lan{\langle}
\def\ran{\rangle}

\def \C{\mathbb{C}}

\def \A{\mathcal{A}}
\def \D{\mathcal{D}}
\def \ov{\overline}

\def \P{\mathbb{P}}
\def \N{\mathcal{N}}
\def \A{\mathcal{A}}

\def \X{X}

\def \lm{\lambda}
\def \G{\gamma}

\theoremstyle{plain}
\newtheorem{thm}{Theorem}[section]
\newtheorem{cor}[thm]{Corollary}
\newtheorem{lmm}[thm]{Lemma}

\newtheorem{que}[thm]{Question}
\newtheorem{eg}[thm]{Example}
\newtheorem{defn}[thm]{Definition}

\newtheorem{rem}[thm]{Remark}

\newtheorem*{ack}{Acknowledgements}

\begin{document}

\maketitle

\begin{abstract}
In this paper we obtain an explicit formula for the number of hypersurfaces 
in a compact 
complex manifold $\X$ 
(passing through the right number of points),
that has a simple node, a cusp or a tacnode. The hypersurfaces 
belong to a linear system, which is obtained by considering 
a holomorphic line bundle $L$ over $\X$. 
Our main tool is a classical fact from 
differential topology: the number of zeros of a generic 
smooth section of a vector bundle $V$ over $M$, counted with a sign, is 
the Euler class of $V$ evaluated on the fundamental class of $M$. 
\end{abstract}

\tableofcontents

\section{Introduction}
\hf\hf Enumeration of singular curves in $\P^2$ (complex projective space) is a classical 
problem in algebraic geometry.
For certain singularities $\mathfrak{X}$, 
the authors 
in \cite{Z1} and \cite{BM13} use a purely 
\textit{topological} method to obtain an explicit  answer for the following question: 
\begin{que}
\label{eg_sing_s}
Let $\mathfrak{X}$ be a codimension $k$-singularity.\footnote{By codimension we mean the number of 
conditions having that particular singularity imposes on the space of curves. 
For example a node is a codimension one 
singularity, a cusp is a codimesnion two singularity, a triple point is a codimension four 
singularity and so on.} 
How many degree $d$-curves are there in $\P^2$, passing through $(d(d+3)/2 -k)$ generic points and having  
a singularity of type $\mathfrak{X}$?
\end{que} 

\ni However, one of the power of that topological method is that it generalizes to 
enumerating curves on \textit{any} complex surface. 
More generally, 
the topological method generalizes 
to enumerating singular hypersurfaces on an arbitrary compact 
complex manifold of a given dimension. \\ 

\ni Let us first make a definition.
\begin{defn}
\label{singularity_defn}
Let $L\lra X$ be a holomorphic line bundle over an $m$-dimensional 
complex manifold $X$ and 
$f:\X \lra L$ a holomorphic section. 
A point $q \in f^{-1}(0)$ \textsf{is of singularity type} $\A_k$ 
if there exists a coordinate system
$(z_1,z_2, \ldots, z_m) :(\mathcal{U},q) \lra (\C^m,0)$ such that $f^{-1}(0) \cap \mathcal{U}$ is given by 
\begin{align*}
z_1^{k+1} + \sum_{i=2}^{m} z_i^2 =0,
\end{align*}
where $m$ is the complex dimension of the manifold $X$.
\end{defn}
\ni Our main theorem can be summarized as follows:
\begin{thm}
Let $\X$ be an $m$-dimensional 
compact complex manifold and $L \lra \X$ a holomorphic line bundle. 
Let 
\[ \mathcal{D} :=\P H^0(\X,L) \approx \P^{\delta_L}, \qquad c_1:= c_1(L) \qquad \textnormal{and} \qquad  x_i:= c_i(T^*\X) \] 
where $c_i$ denotes the $i^{th}$ Chern class. 
Denote $\N(\mathfrak{X})$ to be the 
number of hypersurfaces in $\X$, that belong to the linear system 
$H^0(\X, L)$, passing through $\delta_L-k$ generic 
points and having a singularity of type $\mathfrak{X}$, where $k$ is the codimension of the 
singularity $\mathfrak{X}$. Then 
\begin{align}
\N(\A_1) 
         &= \sum_{i=0}^{m} (m+1-i) x_i c_1^{m-i}, \label{na1_S} \\
\N(\A_2) &= \sum_{i=0}^m m \binom{m+2-i}{2} x_i c_1^{m-i} + \sum_{i=0}^{m-1} 2 \binom{m+1-i}{2} x_1 x_i c_1^{m-i-1}, 
\label{na2_S} \\
\N(\A_3) &= \sum_{i=0}^{m-2} t_2 c_1^{m-i-2} x_i + 
\sum_{i=0}^{m-1} t_1 c_1^{m-i-1} x_i + \sum_{i=0}^{m} t_0 c_1^{m-i} x_i,   \label{na3_S} \\
\textnormal{where} \qquad t_2 &:= \binom{m+1-i}{3} \Big( \frac{c_1^2}{2} (3m^2 -m) + c_1 x_1 (6m-1)  
      + 6 x_1^2 \Big), \nonumber \\
      t_1 &:= \binom{m+1-i}{2} \Big(c_1 (3m^2-m) + x_1 (6m-1) \Big) \qquad \textnormal{and} \nonumber \\
      t_0 &:=\binom{m+1-i}{1} \Big(\frac{3 m^2}{2} -\frac{m}{2} \Big) \nonumber 
\end{align}
provided the sections $\psi_{\A_0}$, $\psi_{\A_1}$, $\psi_{\A_2}$ and 
$\Psi_{\mathcal{P}\A_3}$ 
defined in  
\eqref{psi_a0}, \eqref{psi_a1}, \eqref{psi_a2} and \eqref{psi_pa3} 
are transverse to the zero set, respectively. 
\end{thm}

\begin{rem}
In equations \eqref{na1_S} to \eqref{na3_S} we are making an obvious abuse of notation;
the lhs is an integer, while the rhs is a cohomology class. 
Our intended meaning is the rhs evaluated on the fundamental 
class $[\X]$. 
\end{rem}


\section{Topological computations} 
The main tool that we will use 
is the following well known fact from differential topology 
(cf. \cite{BoTu}, Proposition 12.8).
\begin{thm} 
\label{Main_Theorem} 
Let $V\lra X$ be an oriented  vector bundle over a compact manifold $X$ 
and $s:X \lra V$ a smooth section that is transverse to the zero set. 
Then 
the zero set of $s$ defines an integer homology class  in $X$, 
whose Poincar\'{e} dual is the Euler class of $V$. 
In particular, if the rank of $V$ is same as the dimension of $X$, 
then the signed cardinality of $s^{-1}(0)$ is the Euler class of $V$, evaluated on the fundamental class of 
$X$, i.e., 
\bgd
| \pm s^{-1}(0)| = \lan e(V), [X] \ran. 
\edd
\end{thm}
An immediate corollary is:
\begin{cor}
Let $X$ be a compact, complex manifold, $V$ a holomorphic vector 
bundle (with their natural orientations) 
and $s$ a holomorphic section that is transverse to the zero set. If the rank of $V$ is same as the dimension of $X$, 
then the signed cardinality of $s^{-1}(0)$ is same as its actual cardinality. In particular 
\bgd
|s^{-1}(0)| = \lan e(V), [X] \ran. 
\edd
\end{cor}

\ni We are now ready to give a proof of formulas 
\eqref{na1_S}, \eqref{na2_S} and  \eqref{na3_S}.  
Let $\D_n \approx \P^n \subset \D$ 
be the space of hypersurfaces through 
$\delta_L-n$ generic points. Then $\N(\A_1)$ is the cardinality of the set 
\begin{align}
\{ ([f], q) \in \D_1 \times \X: f(q) =0, ~~\nabla f|_{q} =0\}. 
\end{align}
Let us now define sections of the following bundles: 
\begin{align}
\psi_{\A_0}: \D \times \X \lra &\mathcal{L}_{\A_0}:= \gamma_{\D}^* \otimes L, \qquad \{\psi_{\A_0}([f], q)\}(f) := f(q), 
\label{psi_a0} \\
\psi_{\A_1}: \psi_{\A_0}^{-1}(0) \lra & \mathcal{V}_{\A_1} := \gamma_{\D}^* \otimes T^*\X \otimes L, \qquad \{\psi_{\A_1}([f], q)\}(f) := 
\nabla f|_{q} \label{psi_a1}.
\end{align}
Here $\gamma_{\D}$ is the tautological line bundle over $\D$. 
In equation \eqref{psi_a0}, the rhs is an element of the vector space $L_q$, the fiber of the 
line bundle $L$ over $q$. Hence $\psi_{\A_0}$ is a section of $\mathcal{L}_{\A_0}$. 
Similarly, the rhs of \eqref{psi_a1} is an element of $T^*_q\X \otimes L_q$. Hence  
$\psi_{\A_1}$ is a section of $\mathcal{V}_{\A_1}$.\\
\hf \hf Next, let us assume that $\psi_{\A_0}$ and  $\psi_{\A_1}$ 
are transverse to the zero set.\footnote{Whether or not this assumption actually holds, will depend on 
the specific example of $L$ and $\X$.} 
Since these two sections are holomorphic, we conclude that 
\begin{align}
\N(\A_1) &= \lan e(\mathcal{V}_{\A_1}), ~[\psi_{\A_0}^{-1}(0)] \ran 
         = \lan e(\mathcal{L}_{\A_0})e(\mathcal{V}_{\A_1}), ~[\D_1 \times \X] \ran. \label{na1_euler}
\end{align}
The second equality follows from the fact that the Poincar\'{e} Dual of $[\psi_{\A_0}^{-1}(0)]$ 
in $\D_1 \times \X$ 
is 
$e(\mathcal{L}_{\A_0})$. 
Using the splitting principal, \eqref{na1_euler} simplifies to \eqref{na1_S}. \qed \\
%

\hf \hf Next, let us prove \eqref{na2_S}. 
Note that $\N(\A_2)$ is the cardinality of the following set 
\begin{align}
\{ ([f], q) \in \D_2 \times \X: f(q) =0, ~~\nabla f|_{q} =0, ~~\textnormal{det} \nabla^2 f|_q =0\}. 
\end{align}
Continuing with the setup of the proof of \eqref{na1_S},  
we define a section of the following line bundle 
\begin{align}
\psi_{\A_2}: \psi_{\A_1}^{-1}(0) \lra \mathcal{L}_{\A_2} & :=  
\gamma_{\D}^{*m} \otimes (\Lambda^m T^* \X)^{\otimes 2} \otimes L^{\otimes m},  
\end{align}
given by 
\begin{align}
\label{psi_a2}
\{ \psi_{\A_2}([f], q) \}(f^{\otimes m} \otimes (v_{i_1}\wedge 
\ldots v_{i_m})^{\otimes 2}) & := \textnormal{det} \left( \begin{array}{ccc}
                                          \nabla^2 f|_q(v_{i_1}, v_{i_1})   
                                          & \ldots & \nabla^2 f|_q(v_{i_1}, v_{i_m}) \\ 
                                          \nabla^2 f|_q(v_{i_2}, v_{i_1})   
                                          & \ldots & \nabla^2 f|_q(v_{i_2}, v_{i_m})\\
                                          \ldots  & \ldots & \ldots \\
                                          \nabla^2 f|_q(v_{i_m}, v_{i_1})   
                                          & \ldots & \nabla^2 f|_q(v_{i_m}, v_{i_m})
                                         \end{array}\right).
\end{align}
Note that in \eqref{psi_a2}, the rhs is an element of $L_q^{\otimes m}$. Hence 
$\psi_{\A_2}$ is a section of $\mathcal{L}_{\A_2}$. 
Assume that this section is transverse to the zero set. 
Since it  is holomorphic, we conclude that  
\begin{align}
\N(\A_2) &=  \lan e(\mathcal{L}_{\A_2}), ~~[\psi_{\A_1}^{-1}(0)] \ran \nonumber \\ 
         &= \lan e(\mathcal{L}_{\A_0}) e(\mathcal{V}_{\A_1})e(\mathcal{L}_{\A_2}), ~~[\D_2 \times S]\ran. 
         \label{na2_euler}
\end{align}
The second equality follows from the fact that the Poincar\'{e} Dual of $[\psi_{\A_1}^{-1}(0)]$ 
in $\D_2 \times \X$ 
is $e(\mathcal{L}_{\A_0})e(\mathcal{V}_{\A_1})$. 
Using the splitting principal, \eqref{na2_euler} simplifies to \eqref{na2_S}. \qed \\

\hf \hf Next, we will give a proof of \eqref{na3_S}. 
This requires a bit more ingenuity. First, let us compute $\N(\A_2)$ in an alternate 
way. Let $\P T\X$ denote the projectivization of $T\X$ and  $\hat{\gamma} \lra \P T\X$ 
the tautological line bundle over $\P T\X$. Then $\N(\A_2)$ is also the cardinality 
of the set 
\begin{align}
\{ ([f], l_{q}) \in \D_2 \times \P T\X: (f,q) \in \psi_{\A_1}^{-1}(0), ~~\nabla^2 f |_{q}(v, \cdot) =0 ~~\forall 
~v \in l_{q}\}. \label{psi_pa3} 
\end{align}
Let 
$\pi: \D \times \P T\X \lra \D \times \X$
be the projection map. 
We now define a section of the following bundle 
\begin{align}
\Psi_{\mathcal{P} \A_2}:  \pi^* \psi_{\A_1}^{-1}(0) \lra \mathbb{V}_{\mathcal{P} \A_2} & := 
\hat{\gamma}^* \otimes \gamma_{\D}^* \otimes T^*\X \otimes L, \qquad \textnormal{given by} \nonumber \\ 
\{\Psi_{\mathcal{P} \A_2}([f], l_{q})\}(v \otimes f) & := \nabla^2 f|_{q}(v, \cdot) \qquad \forall ~~v \in l_q.
\label{psi_a2_section_defn}
\end{align}
Note that the rhs of \eqref{psi_a2_section_defn} belongs to $T^*_q\X \otimes L_q$. Hence 
$\Psi_{\mathcal{P}\A_2}$ is a section of $\mathbb{V}_{\mathcal{P} \A_2}$. 
Assume that this section 
is transverse to the zero set. 
Since it is holomorphic, we conclude that 
\begin{align}
\N(\A_2) &= \lan e(\mathbb{V}_{\mathcal{P}\A_2}), ~[\pi^* \psi_{\A_1}^{-1}(0)] \ran \nonumber \\
         &= \lan e(\pi^* \mathcal{L}_{\A_0}) e( \pi^*\mathcal{V}_{\A_1})
         e(\mathbb{V}_{\mathcal{P}\A_2}), ~[\D_2 \times \P T\X] \ran.  \label{n_pa2_euler}
\end{align}
We now use the ring structure of $H^*(\P T\X, \mathbb{Z})$  
(\cite{BoTu}, pp. 270) to conclude that  
\begin{align}
\lambda ^m - c_1 \lambda^{m-1} + c_2 \lambda^{m-2} + \ldots + (-1)^m c_m &=0, \label{lambda_square} 
\end{align}
where $\lambda := c_{1}(\hat{\gamma}^*)$.
Using the splitting principal, equations \eqref{n_pa2_euler} combined with  \eqref{lambda_square} simplifies 
to \eqref{na2_S}. \\

\hf \hf We are now ready to prove \eqref{na3_S}.   
First, we note that $\N(\A_3)$ is also the cardinality of the set 
\begin{align*}
\{ ([f], l_{q}) \in \D_3 \times \P T\X: ([f], l_q) \in 
\Psi_{\mathcal{P}\A_2}^{-1}(0), ~~\nabla^3 f |_{q}(v, v,v) =0 ~~\forall 
~v \in l_{q}\}. 
\end{align*}
Let us now define a section of the following bundle 
\begin{align}
\Psi_{\mathcal{P} \A_3}:  \Psi_{\mathcal{P}\A_2}^{-1}(0) 
\lra \mathbb{L}_{\mathcal{P} \A_3} & := \hat{\gamma}^{*3} \otimes \gamma_{\D}^* \otimes L, 
\qquad \textnormal{given by} \nonumber \\
\{\Psi_{\mathcal{P} \A_3}([f], l_{q})\}(v^{\otimes 3} \otimes f) & := \nabla^3 f|_{q}(v, v,v) \qquad 
\forall ~~v \in l_q. \label{psi_pa3_section_defn}
\end{align}
Note that 
the rhs of \eqref{psi_pa3_section_defn} is an element of $L_q$. Hence 
$\Psi_{\mathcal{P} \A_3}$ is a section of $\mathbb{L}_{\mathcal{P} \A_3}$.
Assume that 
this section is transverse to the zero set. 
Since it is holomorphic, we conclude that 
\begin{align}
\N(\A_3) &= \lan e(\mathbb{L}_{\mathcal{P} \A_3}), ~[\Psi_{\mathcal{P}\A_2}^{-1}(0)] \ran \nonumber  \\ 
         &= \lan 
         e(\pi^* \mathcal{L}_{\A_0}) e( \pi^*\mathcal{V}_{\A_1})
         e(\mathbb{V}_{\mathcal{P}\A_2})e(\mathbb{L}_{\mathcal{P} \A_3}), ~[\D_3 \times \P TS] \ran. 
         \label{n_pa3_euler} 
\end{align}
Using the splitting principal, equations \eqref{n_pa3_euler} combined with  \eqref{lambda_square} simplifies 
to \eqref{na3_S}. \qed

\section{Examples}
\begin{eg}
\label{na1_surface_var} 
Let $\X$ be a surface. Then $\eqref{na1_S}$ simplifies to 
\begin{align}
\N(\A_1) &= 3 c_1^2 + 2 c_1 x_1 + x_2, 
\end{align}
which agrees with the result of Vainsencher in \cite{Van}.
\end{eg}

\begin{eg}
\label{eg_pm}
Let $\X:= \P^m$ and $L:= \gamma_{\P^m}^{*d}$. Then  
\begin{align}
c_1^{m-i} x_i &= (-1)^i\binom{m+1}{i} d^{m-i}. \label{c1_xi_na1_pm}
\end{align}
Equations \eqref{na1_S}, \eqref{na2_S} and \eqref{na3_S} 
combined with \eqref{c1_xi_na1_pm}   imply that 
\begin{align}
\N(\A_1) & = (m+1) (d-1)^m,  \\
\N(\A_2)  &= \frac{m(m+1)(m+2)}{2} (d-1)^{m-1} (d-2), \\
\N(\A_3) &= \frac{m(m+1)(m+2)}{12} (d-1)^{m-2} (m_2 d^2 + m_1 d + m_0),\\  
\textnormal{where} \qquad m_2 &:= m^2+2m-1, \qquad m_1:= -12m^2-28m+8 \qquad \textnormal{and} \qquad m_0:= 3m^2 + 8m-3. \nonumber 
\end{align}
For $m=2$ this gives us  
\begin{align}
\N(\A_1) &= 3(d-1)^2, ~~\N(\A_2) = 12(d-1)(d-2), ~~\N(\A_3) = 50d^2-192d+168, 
\end{align}
which recovers the formulas obtained in \cite{Z1} and \cite{BM13}. 
For general $m$, 
the numbers $\N(\A_2)$ and $\N(\A_3)$  agree with the results of Kerner in 
\cite{Ker_Hypersurface}.\footnote{The number $\N(\A_1)$ is not explicitly stated in \cite{Ker_Hypersurface}.} 
\end{eg}


\begin{eg}
\label{eg_p1_times_p1}
Let $\X:= \P^1 \times \P^1$ and $L:= \pi_1^* \gamma^{* d_1}_{\P^1} \otimes \pi_2^* \gamma^{* d_2}_{\P^1}$.  
Then 
\begin{align}
c_1^2 &=2 d_1 d_2 , ~~c_1 x_1 = -2(d_1 +d_2) , 
~~x_1^2 = 8 , ~~x_2 = 4. \label{chern_num_P1_times_p1}
\end{align}
Equation \eqref{chern_num_P1_times_p1}, combined with \eqref{na1_S}, \eqref{na2_S} and \eqref{na3_S} gives 
\begin{align}
\N(\A_1) &= 6d_1 d_2 - 4(d_1+d_2) +4 , \nonumber \\ 
\N(\A_2) & = 24(d_1-1)(d_2-1) \qquad \textnormal{and}    \nonumber  \\ 
\N(\A_3) &=  100 d_1 d_2 -128 (d_1 + d_2) + 136.  \nonumber  \\
\end{align}
Let us take $d_1 = d_2 =1$. The formulas for $\N(\A_2)$ and $\N(\A_3)$ confirm the fact that there 
are no curves of type $(1,1)$ with either a cusp or a tacnode. The formula for $\N(\A_1)$ confirms the 
fact that there are two curves of type $(1,1)$ through two generic points that have a node. 
\end{eg}

\begin{rem}
In \cite{BM_Detail}, we show why the relevant sections in example \ref{eg_pm} 
are transverse to the 
zero set (for $m=2$). 
The proof follows by unwinding definitions and writing the section in a local coordinate 
chart and trivialization. The proof of why the relevant sections  
in example \ref{eg_p1_times_p1} are transverse 
to the zero set is similar (provided $d_1$ and $d_2$ are sufficiently large). 
\end{rem}

\begin{ack}
The author thanks Indranil Biswas for suggesting example \ref{eg_p1_times_p1}. 
Furthermore, the author is also grateful to Somnath Basu, Shane D'Mello and Vamsi Pingali 
for relevant discussions and comments about this paper. 
\end{ack}

\bibliographystyle{siam}
\bibliography{Myref_bib.bib}

\begin{thebibliography}{1}

\bibitem{BM13}
{\sc S.~Basu and R.~Mukherjee}, {\em Enumeration of curves with one singular
  point}.
\newblock available at \url{http://arxiv.org/abs/1308.2902}.

\bibitem{BM_Detail}
\leavevmode\vrule height 2pt depth -1.6pt width 23pt, {\em Enumeration of
  curves with singularities: Further details}.
\newblock available at \url{https://www.sites.google.com/site/ritwik371/home}.

\bibitem{BoTu}
{\sc R.~Bott and L.~W. Tu}, {\em Differential forms in algebraic topology},
  vol.~82 of Graduate Texts in Mathematics, Springer-Verlag, New York, 1982.

\bibitem{Ker_Hypersurface}
{\sc D.~Kerner}, {\em Enumeration of uni-singular algebraic hypersurfaces},
  Proc of London Mathematical Society, 3 (2008), pp.~623--668.

\bibitem{Van}
{\sc I.~Vainsencher}, {\em Enumeration of {$n$}-fold tangent hyperplanes to a
  surface}, J. Algebraic Geom., 4 (1995), pp.~503--526.

\bibitem{Z1}
{\sc A.~Zinger}, {\em Counting plane rational curves: old and new approaches}.
\newblock available at \url{http://arxiv.org/abs/math/0507105}.

\end{thebibliography}


\begin{thebibliography}{99}

\bibitem{BM1} S.Basu, R.Mukherjee,   ~\textit{Enumeration of curves with one singular point},  
available at \url{https://www.sites.google.com/site/ritwik371/home}.


\bibitem{BM_Detail} S.Basu, R.Mukherjee,   ~\textit{Enumeration of curves: further details},  
available at \url{https://www.sites.google.com/site/ritwik371/home}.


\bibitem{Z1} A.Zinger,   ~\textit{Counting plane rational curves: old and new approaches},  
available at \url{http://arxiv.org/abs/math/0507105}.

\end{thebibliography}

\vspace*{0.4cm}


\hf {\small D}{\scriptsize EPARTMENT OF }{\small M}{\scriptsize ATHEMATICS, }{\small TIFR,}
{\small T}{\scriptsize ATA  }
{\small I}{\scriptsize NSTITUTE }
{\small OF}
{\scriptsize FUNDAMENTAL}{\small R} {\scriptsize ESEARCH,} 
{\small M}{\scriptsize UMBAI }
{\footnotesize 400005, }{\small INDIA}\\
\hf{\it E-mail address} : \texttt{ritwikm@math.tifr.res.in}\\[0.2cm]

\end{document}